\magnification=\magstep1
\def\ds{\baselineskip 18pt plus 2pt}
\def\ss{\baselineskip 10pt plus 1pt}

\def\chs {compact Hausdorff space \enspace}
\def\bs {Banach space \enspace}
\def\bss {Banach spaces \enspace}
\def\as {absolutely summing \enspace}
\def\aso {absolutely summing operator\enspace}

\def\bv {bounded variation\enspace}

\def\io { integral operator\enspace }
\def\ios { integral operators\enspace }
\def\blo { bounded linear operator\enspace }
\def\blos { bounded linear operators\enspace }

\def\intnorm#1{\parallel#1\parallel_{\text {int}}}

\def\etensorf { E \hat \otimes_\epsilon F}
\def\etensorfstar { E \hat \otimes_\epsilon F^*}
\def\emu{\langle e,\mu \rangle}
\def\mxestar {M(X,E^*)}

\def\onebe {1_B \otimes e}
\def\cxe {C(X,E)}
\def\cxbestar {C(X \times B(E^*))}
\def\cbestar {C( B(E^*))}
\def\xbestar {X \times B(E^*)}

\def\gystar{G_{y^*}}
\def\cbestarbfdoublestar{C(B(E^*)\times B(F^{**}))}
\def\cxbestarbfdoublestar{C(X\times B(E^*)\times B(F^{**}))}
\def\xbestarbfdoublestar{X\times B(E^*)\times B(F^{**}}
\def\lefdoublestar{{\Cal L}(E,F^{**})}
\def\bestarfdoublestar{B(E^*)\times B(F^{**})}
\def\bestarbfdoublestar{B(E^*)\times B(F^{**})}
\def\bestar { B(E^*)}
\def\linfinitymu {L_{\infty}(\mu)}
\def\lonemu {L_{1}(\mu)}

\def\mx {M(X)}
\def\cx {C(X)}
\def\norm#1{\|#1\|}
\def\abs#1{\vert #1\vert}  
\def\sequence #1#2#3{(#1_#2)_{#2\ge #3}} 
\def \theorem #1#2{\medbreak\noindent{\bf#1\enspace}{\sl#2}\medbreak
\ifdim\lastskip<\medskipamount \removelastskip\penalty55\medskip\fi}
\def\equal{\,=\,}
\def\inn {\, \in \,}
\def\map #1#2#3{#1 \colon #2 \longrightarrow #3}

\pageno=0
\footline={\ifnum\pageno=0\hfill\else\hss\tenrm\folio\hss\fi}

\font\bigbold=cmbx10 scaled \magstep2
\def\ans{\vrule height.1pt width80pt depth0pt}

\voffset=1in
\centerline {\bigbold Integral Operators on Spaces of}
\centerline {\bigbold Continuous Vector-valued functions}
\vskip 1truein
\centerline {by}
\vskip .50truein
\centerline {\bf Paulette Saab$^*$}
\vskip 1truein
\ss

{\narrower\smallskip\noindent {\bf Abstract}\ \ Let $X$ be a compact
Hausdorff space, let
$E$ be a Banach space, and let $C(X,E)$ stand for the Banach space of
$E$-valued continuous functions on $X$ under the uniform norm.
In this paper we
characterize Integral operators (in the sense of Grothendieck) on
$C(X,E)$ spaces in term of their representing vector measures.  This is
then used to give some applications to Nuclear operators on $C(X,E)$
spaces.\smallskip}
\vskip 2truein
{\narrower\smallskip\noindent AMS(MOS) subject Classification (1980).
Primary
46E40,
46G10;\ Secondary 28B05, 28B20.\smallskip}

\ans

$*$ Supported in part by an NSF Grant DMS-87500750
\vfill\eject
\voffset=-.3in
\ds

\noindent {\bf Introduction}\ \ Let $X$ be a compact Hausdorff space,
let $E$ and $F$ be Banach spaces. Denote by $\cxe$ the space of all
continuous $E$-valued functions defined on $X$ under the uniform norm. In
[9] C. Swartz showed that a bounded linear operator
$\map T {\cxe} F$ with representing measure $G$ is absolutely summing if
and only if each of the values of $G$ is an \aso from $E$ to $F$ and $G$
is of \bv  as a measure taking its values in the space of \as operators
from $E$ to $F$ equipped with the \as norm. In this paper we shall extend
Swartz's result to the class of (Grothendieck) \ios on $\cxe$ spaces.
More precisely we shall show that a \blo
$T:\ C(X,E)\rightarrow F$ with representing measure $G$
is an integral operator if and only if each of  the values of $G$ is an
\io from $E$ to $F$ and $G$
is of bounded variation as a vector measure
taking its values in the space of \ios from $E$ to $F$ equipped
with the integral norm.
This result is then  used to give some applications to Nuclear operators
on
$C(X,E)$ spaces.
\bigskip
\noindent {\bf I. Preliminaries}\  If $X$ is a \chs and $E$ is a \bs, then
$\cxe$ will denote the \bs of all continuous $E$-valued functions
equipped with the uniform norm. It is well known [4,page 182] that the
dual of $\cxe$ is isometrically isomorphic to the space $\mxestar$
of all regular $E^*$-valued measures on $X$ that are of bounded
variation. When $E$ is the scalar field, we will simply write $\cx$ and
$\mx$ for $\cxe$ and $\mxestar$. If $\mu \inn \mxestar$ and $e \inn E$,
we
will denote by $|\mu|$ the variation of $\mu$ and by $\emu$
the element of $\mx$ defined on each Borel subset $B$ of $X$ by
$$
\emu (B) \equal \mu (B) (e)
.$$
The duality between $\mxestar$ and $\cxe$ is then defined as follows: for
each $f\inn \cx$ and $e \inn E$
$$
\mu (f \otimes e) \equal \int_X f d \emu
$$
where $ f \otimes e $ is the element of $\cxe$ defined by
$$
f\otimes e (x) \equal f(x) e \text { for all } x \inn X .
$$

If $B$ is a Borel subset of $X$, then $1_B$ will denote the
characteristic function of $B$, and if $e\inn E$ we let
$\onebe$ denote the element of $\cxe^{**}$ defined by
$$
\onebe (\mu) \equal \emu (B) \equal \mu (B)(e)
$$
for each $\mu \inn \mxestar$.

If $X$ is a \chs , $E$ and $F$ are \bss, every \blo
$\map T {\cxe} F$ has a representing measure $G$. The measure $G$ is
defined on the $\sigma$-field $\Sigma$ of Borel subsets of $X$ and takes
its values in
$\lefdoublestar$, the space of all \blos from $E$ to $F^{**}$. The
measure $G$ is such that for each Borel subset $B$ of $X$ and for each
$e\inn E$
$$
G(B) e \equal T^{**} (\onebe)
$$
For $y^* \inn F^*$, if we denote by $\gystar$ the $E^*$- valued measure on
$X$ such that for each Borel subset $B$ of $X$ and each $e \inn E$
$$
\langle e,\gystar \rangle (B) \equal \langle y^*, G(B) e \rangle
$$
then $\gystar$ is the unique element of $\mxestar$ that represents $T^*
y^*$ in the sense that for each $f\inn \cxe$
$$
\langle y^*, Tf \rangle \equal \int_X f(x) d\gystar (x)
$$
If $E$ ad $F$ are Banach space, we denote by $E \otimes_\epsilon F$ the
algebraic tensor product of $E$ and $F$ endowed with the norm
$|| .||_{\epsilon }$
$$
||\sum_{i=1}^m x_i \otimes y_i ||_{\epsilon} \equal \sup \{|\sum_{i=1}^m
x^*(x_i) y^*(y_i) | \quad \vert \quad || x^* ||\, , \, || y^* || \leq 1\}.
$$
The completion $\etensorf$ of $E\otimes_\epsilon F$ is called the
injective tensor product of $E$ and $F$.

A \blo
$T:\ E\rightarrow F$ from a Banach space
$E$ into a Banach space $F$ is said to be an {\bf integral operator} if
the bilinear form $\tau$ on $E\times F^*$ defined by $\tau
(e,y^*)=y^*(Te)$ for $e\in E$ and $y^*\in F^*$ determines an element of
$(E\hat \otimes_\epsilon F^*)^*,$ the dual of the injective tensor
product
$E\hat \otimes_\epsilon F^*$
of the Banach spaces $E$ and $F^*$.  The integral norm of $T$, which we
will denote by $\parallel T \parallel_{\text {int}}$, is just the
norm of the bounded linear functional induced by $\tau$ as an element of
the dual space $(E\hat \otimes_\epsilon F^*)^*$.  Hence if $T:\
E\longrightarrow F$ is an integral operator, then
$$\parallel T\parallel_{\text {int
}}=\sup\left\{|\sum\limits^n_{i=1}y^*_i(T e_i)|:\ \parallel
\sum\limits^n_{i=1}e_i\otimes y_i^*\parallel_\epsilon \leq 1\right\}
.$$
Finally,
$I(E,F)$ will stand for
the \bs of all integral operators from $E$ to $F$ equipped with the
integral norm. For
all undefined notions and notations we refer the reader to [3], [4] or
[5].
\bigskip

\ One of the most useful tools in the
study of an integral operators $T$ between two Banach spaces $E$ and $F$
is its representation by a finite regular Borel scalar measure on the
compact $B(E^*)\times B(F^{**})$, the product of the
closed unit balls of $E^*$
and $F^{**}$ equipped with their weak$^*$-topologies.  This, of course,
goes back to Grothendieck [6] and can be used to characterize integral
operators by the following characterization that we will state and prove
before proving the main result.

\bigskip
\theorem { Proposition 1} { A bounded linear operator $T:\
E\rightarrow F$ between two Banach spaces is an integral operator if and
only if there exists a regular $F^{**}$-valued vector measure $m$ of \bv
defined
on the $\sigma$-field of  Borel subsets of the closed unit
ball $B(E^*)$,
such that for each $e\in E$
$$
Te=\int_{B(E^*)}e^*(e)\ dm(e^*)
$$
In this case $m$ can be chosen so that $\intnorm T \equal \abs m
(\bestar).$ }
\bigskip
\noindent
{\bf Proof:} Assume that there exists a regular $F^{**}$-valued vector
measure $m$ defined on the $\sigma$-field of  Borel subsets of
$\bestar$ with
$\abs m (\bestar) < \infty$ and such that for each $e\inn E$
$$
Te=\int_{B(E^*)}e^*(e)\ dm(e^*)
$$
This in particular shows that the operator
$\map {J\circ T} E {F^{**}}$, where $J$ denotes the natural embedding of
$F$ into $F^{**}$ , has an extension $\hat T$ to an \io from $\cbestar$
to $F^{**}$. Hence $T$ is an \io [4,page 233] and
$$
\intnorm T \leq \intnorm {\hat T} \equal \abs m (\bestar)
.$$

Conversely, suppose that $\map TEF$ is an \io, then it follows from
[4,page~231] that there exists a regular Borel measure $\mu$ on
$\bestarbfdoublestar$ such that for each $e\inn E$ and $y^* \inn F^*$
$$
\langle y^*, Te \rangle \equal
\int _{\bestarbfdoublestar} e^*(e) y^{**} (y^*) d\mu(e^*,y^{**})
$$
and
$$
\intnorm T \equal \abs \mu (\bestarbfdoublestar)
.$$
Following [4,page 234], define
$$\map S E {\linfinitymu} \quad\text { by }\quad Se (e^*,y^{**}) \equal
e^*(e)$$ and
$$\map R {F^*} {\linfinitymu}\quad \text { by }\quad R y^* (e^*,y^{**})
\equal y^{**}(y^*).$$
Then $S$ and $R$ are \blos with
$\norm S, \norm R \leq 1$. Let $Q$ be the restriction of $R^*$ to
$\lonemu$, then $Q$ is a \blo from $\lonemu$ into $F^{**}$. It is
immediate that $J\circ T \equal Q\circ I \circ S$ where
$\map I {\linfinitymu} {\lonemu}$ is the natural inclusion
and $\map J F {F^{**}}$ is the natural embedding. This in particular
shows that the operator $Q\circ I$ restricted to $\cbestarbfdoublestar$
is an \io whose representing $F^{**}$-valued measure $\hat \mu$ is such
that
$$
\hat \mu (C) \equal Q(1_C)
$$
for each Borel subset $C$ of $\bestarfdoublestar$. Finally note that
since the natural projection of $\bestarfdoublestar$ onto $\bestar$ is a
continuous mapping it induces a \blo from $\cbestar$ into
$\cbestarbfdoublestar$ as follows; for $\phi \inn \cbestar$, let $\hat
\phi \inn \cbestarbfdoublestar$ be such that
$$
\hat \phi ( e^*,y^{**}) \equal \phi (e^*)
.$$
This of course shows that the operator $J \circ T$ extends to an integral
operator
$$\map {\hat T} {\cbestar} {F^{**}}$$ such that for each $\phi \inn
\cbestar$
$$
\hat T \phi \equal Q \circ I (\hat \phi)
$$
It is immediate that the $F^{**}$-valued measure $m$ representing $\hat
T$ is such that
$$
m(B) \equal Q(1_{B \times B(F^{**})})
$$
and
$$
\abs m (\bestar \leq \abs \mu (\bestarfdoublestar) \equal \intnorm T
.$$
Since $\hat T$ extends $T$ it follows that
$$
\intnorm T \equal \abs m (\bestar).
$$

\noindent{\bf Main Result} Throughout this section $X$ is a compact
Hausdorff space,
$E$ and
$F$ are \bss and $\map T {\cxe} F$ is a \blo with representing measure
$G$. The main result of this paper gives a characterization of \ios $T$
in terms of some properties of $G$. The first step to achieve such a
characterization is to show that if $T$ is an \io on $\cxe$, then one can
do a little better than Proposition 1 by representing the operator $T$ by
a regular $F^{**}$-valued measure of \bv defined on the $\sigma$-field of
Borel subsets of $\xbestar$ rather than on the whole unit ball of
$({\cxe})^*$. The proof we present here is different from our earlier
proof which relied on our result [7]. We would like to thank the referee
for suggesting the following approach which relies  more on basic
knowledge
and classical results of vector measures that can be found in [4] or
[5].

\theorem { Lemma 2 .} {If $\map T {\cxe} F$ is an \io, then there
exists
a regular $F^{**}$-valued measure of \bv $\theta$ on the $\sigma$-field
of Borel subsets of $\xbestar$ such that for each $f \inn \cxe $ and
$y^* \inn F^*$
$$
\langle y^*, T f \rangle \equal \int_{\xbestar} e^* ( f(x)) d \theta
(x,e^*)
$$
and
$$
\intnorm T \equal \abs \theta (\xbestar).
$$ }
\bigskip\noindent
{\bf Proof:} Suppose that $\map T {\cxe} F$ is an \io, then the bilinear
map $\tau (f,y^*) \equal y^*(T(f)) $ for $f\inn \cxe$ and $y^* \inn F^*$
defines an element of
$( \cxe \otimes_\epsilon F^*)^* $.
 It is easy to check that
$ \cxe \otimes_\epsilon F^* $ embeds isometrically in
$\cxbestarbfdoublestar$ such that if $f\inn \cxe$ and $y^* \inn F^*$
$$
f \otimes y^* (x,e^*,y^{**}) \equal e^*( f(x)) y^{**}(y^*)
$$
for each $(x,e^*,y^{**}) \inn \Omega \equal \xbestarbfdoublestar$.
Hence by the Hahn-Banach theorem there exists a regular Borel measure
$\mu$ on $\Omega$ so that for each $f \inn \cxe $ and $y^* \inn F^*$
$$
\langle y^*, Tf \rangle \equal \int_\Omega
e^* (f(x)) y^{**}(y^*) d\mu (x,e^*,y^{**}).
$$
and
$$
\intnorm T \equal \abs \mu (\Omega).
$$
The proof now follows the steps of Proposition 1. Indeed, since the
natural projection of $\Omega$ onto $\xbestar$ is continuous it induces a
\blo
$$
\map S {\cxbestar} {\linfinitymu}
$$
such that if $\phi \inn \cxbestar$ then
$$
S(\phi)(x,e^*,y^{**}) \equal \phi (x,e^*)
$$
for all $(x,e^*,y^{**}) \inn \Omega .$

Let $\map R {F^*} {\linfinitymu}$ be defined by
$R(y^*) (x,e^*,y^{**}) \equal y^{**}(y^*)$
and let $Q$ be the restriction of $R^*$ to $\lonemu$. It is
straightforward
to check that the operator $Q \circ I \circ S $ is an \io
 whose representing measure $\theta$ is such
that
$$
\theta (C) \equal Q (1_{C \times B(F^{**})})
$$
for each Borel subset $C$ of $\xbestar$, and
$$
\abs \theta (\xbestar) \leq \abs \mu (\Omega) \equal \intnorm T .
$$
Moreover note that $\cxe$ embeds isometrically in $\cxbestar$ by
$$
f(x,e^*) \equal e^*(f(x))
$$
for each $f \inn \cxe$ and $(x,e^*) \inn \xbestar$. It is easy to see
that for each $f \inn \cxe$
$$
J \circ T (f) \equal Q \circ I \circ S (f)
$$
hence
$$
 T(f)  \equal \int_{\xbestar} e^*(f(x)) d\theta(x,e^*)
$$
and
$$
\intnorm T \leq \intnorm {Q \circ I \circ S} \leq \abs \theta (\xbestar)
.$$
This complete the proof of Lemma 2 .

We are now ready to state and prove the main result of this paper.

\theorem { Theorem 3} {Let $X$ be a compact Hausdorff space, let
$E$ and $F$
be Banach spaces and let $T:\ C(X,E)\rightarrow F$ be a bounded linear
operator with representing measure $G$.  Then $T$ is an integral
operator if and only if for
each Borel subset $B$ of $X$, the operator $\map {G(B)} E F$ is
integral and the measure $G$ is of finite variation as a vector measure
taking its values in $I(E,F)$ equipped with the integral norm.}
\bigskip
\noindent {\bf Proof:}\ \ Assume that $T:\ C(X,E)\rightarrow F$ is an
integral operator.  Let $G$ denote the vector measure representing the
operator $T$, hence for each Borel subset $B$ of $X$ and each $e\in E$
$$
G(B)e=T^{**}(\onebe)
$$
It is clear at this stage, that since $T^{**}$ is also integral [4,
p.236
], then for each Borel subset $B$ of $X$ the operator $\map {G(B)} E F $
 is an integral operator as the composition of $T^{**}$ and the
bounded linear operator $E\longrightarrow C(X,E)^{**}$ which to each $e$
in $E$ associates the element $\onebe$.
In what follows we shall concentrate on estimating the value of
$\intnorm {G(B)}$. For this note that by Lemma 2 there exists a regular
$F^{**}$-valued measure $\theta$ defined on the $\sigma$-field of Borel
subsets of $\xbestar$ such that, for each $f \inn \cxe$
$$
T(f) \equal \int_{\xbestar} e^*(f(x)) d\theta(x,e^*)
$$
and
$$
\intnorm T \equal \abs \theta (\xbestar) .
$$
It follows that for each $y^* \in F^*$ and $f\in C(X,E)$
$$
\langle T^*y^*,f\rangle=\int_{X\times B(E^*)}e^*\left(f(x)\right)
d\theta_{y^*}(x,e^*).
$$
where $\theta_{y^*}$ is the scalar measure on $X\times B(E^*)$ defined
by
$$
\theta_{y^*}(C)=\theta(C)(y^*)
$$
for each Borel subset $C$ of $X\times B(E^*)$.

\noindent
We claim that for each Borel subset $B$ of $X$ and for each $e \inn E$
$$
G(B)e \equal \int_{B \times B(E^*)} e^*(e) d\theta(x,e^*) . \tag *
$$
For this suppose that $y^* \inn F^*$ and $K$ is a compact $G_\delta $
subset of
$X$. Let $\sequence un1$ be a sequence of continuous real valued functions
so that $0 \leq u_n \leq 1  \text { for all } n \geq 1$ and
$u_n$ converges to $1_K$ pointwise. If $e \inn E$, then

$$
\aligned
\langle G(K)e,y^*\rangle &=\langle
e,\gystar \rangle (K)\\
&\\
&=\lim\limits_{n\to \infty} \langle y^*, T(u_n \otimes e )\rangle \\
&\\
&=\lim\limits_{n\to \infty} \int_{\xbestar} u_n(x) e^*(e) d\theta_{y^*}
(x,e^*)\\
&\\
&=\int_{\xbestar} 1_K(x)  e^*(e) d\theta_{y^*}
(x,e^*)\\
&\\
&=\int_{K \times B(E^*)}   e^*(e) d\theta_{y^*}
(x,e^*)\\
\endaligned .$$
Therefore
$$
G(K) e \equal
\int_{K \times B(E^*)}   e^*(e) d\theta (x,e^*)
.$$
Moreover, since for each $y^* \inn F^*$ the $E^*$-valued vector measure
$\gystar$ is regular, it follows that
$$
G(B) e \equal
\int_{B \times B(E^*)}   e^*(e) d\theta (x,e^*)
$$
for all Borel subsets $B$ of $X$. This proves our claim.
If we denote by $p:\
X\times B(E^*)\longrightarrow E^*$ the projection mapping which to each
$(x,e^*)$ in $X\times B(E^*)$ associates $e^*$ in $B(E^*)$, then for
each Borel subset $B$ of $X$, let $\lambda_B$ denote
 the regular $F^{**}$-valued
measure defined on the $\sigma$-field of Borel subsets of $B(E^*)$ as
follows:\ \ for each Borel subset $V$ of $B(E^*)$
$$
\lambda_B(V)=\theta (B\times V).
$$

In other words $\lambda_B=\theta|B\times B(E^*)\circ p^{-1}$ which is the
image measure of the restriction of $\theta$ to $B\times B (E^*)$ by
$p$.  This implies that for each $e\in E$
$$
\int_{B(E^*)}e^*(e)d\lambda_B(e^*)=\int_{B\times B(E^*)}e^*(e)\ d\theta
(x,e^*) . \tag **
$$
Equations $(^{*})$ and $(^{**})$ show that the measure $\lambda_B$ is a
regular
$F^{**}$-valued measure that represents the operator $\map {G(B)} E F$.
Hence by Proposition 1
$$
\intnorm {G(B)} \leq \abs {\lambda_B} (\bestar)
\equal \abs \theta ( B \times \bestar) .
$$
and therefore
$$
(\dagger)\qquad \parallel G(B)\parallel_{\text {int}}\leq |\theta|\
(B\times B(E^*)
$$
Hence the vector measure $G$ representing the operator $T$ takes
its values in $I(E,F)$ and it follows from $(\dagger)$ that
$$
\qquad |G|_{\text {int}}(X)\leq |\theta|\left(X\times B(E^*)\right).
$$
Here of course $|G|_{\text {int}}(X)=\sup \sum\limits_{B_i\in
\pi}\parallel
G(B_i)\parallel_{\text {int}}$ where the sup is taken over all the finite
partitions of $X$ into Borel subsets of $X$.
\medskip
Conversely, suppose $T:C(X,E)\longrightarrow F$ is such that
$$\map G
\Sigma {I(E,F)} \quad \text { and } |G|_{\text {int}}(X)<~\infty
,$$ we
need
to show that $T$ is an integral operator.  For this, note that if $J$
denotes the natural embedding of $F$ into $F^{**}$, then by [4, p. 233]
it is enough to show that $J\circ T:\ C(X,E)\longrightarrow F^{**}$ is
an integral operator.  The mapping $J$ induces a mapping $\hat J:\
I(E,F)\longrightarrow I(E,F^{**})$ defined as follows, for each $U\in
I(E,F)$
$$
\hat J(U)=J\circ U
$$
this in turn induces a vector measure $\hat G:\ \Sigma \longrightarrow
I(E,F^{**})$ such that
$$
\hat G(B)=\hat J\left(G(B)\right)
$$
for all $B\in \sum$.  It is immediate that $\hat G$ is the measure
representing the operator $J\circ T$.  Since $|G|_{\text
{int}}(X)<\infty$, it follows easily that $|\hat G|_{\text
{int }}(X)<\infty$.  Moreover, since $I(E,F^{**})$ is isometric to the
dual space of the injective tensor product $E \hat \otimes_{\epsilon}
F^*[4, p.237]$,
it follows that  $\hat G$ is a vector measure of \bv taking its
values in the dual
space $(\etensorfstar)^*$. Since $|\hat G|_{\text {int}} (X) < \infty$,
it
follows from [4,page 7] that $\hat G$ is strongly additive. Recall that
for each $e\inn E$ and $y^* \inn F^*$, the scalar measure $\langle
e,\gystar \rangle$ is in $\mx$. Hence for
$u \equal \sum\limits_{i=1}^n e_i \otimes y_i^* \inn \etensorfstar$, let
$\langle u,\hat G \rangle $ be the scalar measure defined on $X$ as
follows: for each Borel subset $B$ of $X$

$$
\aligned
\langle u,\hat G\rangle(B) &= \hat G(B) (u)
\\
&\\
&=\sum\limits_{i=1}^n \langle G(B) e_i, y_i^* \rangle \\
&\\
&=\sum\limits_{i=1}^n \langle e_i,G_{y^*_i} \rangle (B)
\\
\endaligned .$$
This implies that for each $u \inn \etensorfstar $ the measure
$\langle u, \hat G \rangle \inn \mx $. This in particular shows that
$\map {\hat G} \Sigma {(\etensorfstar)^*}$ is countably additive with
respect to the $\hbox {weak}^* \,$  topology on
$(\etensorfstar)^*$. Since
$\hat G$
is strongly additive, it follows easily that $G$ is countably additive.
Finally, since $\hat G$ is $weak^* \,$ regular, a glance at [4,page
157] reveals that it is also regular.
It follows that $\hat G$ defines an element of $C(X,
E \hat \otimes_{\epsilon} F^*)^*$, hence by [5, page 269] there exists a
weak$^*\ |\hat G|_{\text {int }}$-integrable function
$$
h:\ \ X\longrightarrow I(E,F^{**})\simeq (E\hat \otimes_{\epsilon}
F^*)^*
$$
such that $\parallel h (x)\parallel_{\text {int}}=1\quad |\hat
G|_{\text {int }}$ a.e.,
$\hat G=h|\hat G|_{\text {int}}$, and for each $f\in C(X,
E\hat \otimes_{\epsilon} F^*)$
$$
\langle f, \hat G\rangle =\int_X\langle f(x), h(x)\rangle \ d|\hat
G|_{\text {int }}(x)
$$
In particular, if $\varphi \in C(X,E)$ and $y^*\in F^*$, let
$f=\varphi \otimes y^*$ be the element of $C(X, E\hat \otimes_\epsilon
F^*)$ defined by
$$
\varphi \otimes y^* (x)=\varphi (x)\otimes y^*\text { for all }x\in X.
$$
then
$$
\langle \varphi \otimes y^*, \hat G\rangle=\int_Xy^*\langle h(x),
\varphi (x)\rangle d|\hat G|_{\text {int }}| (x)
$$
But we also have
$$
\langle J\circ T\varphi, y^*\rangle=\langle \varphi, \hat
G_{y^*}\rangle\tag ***
$$
where $\hat G_{y^*}$ is the element of $M(X,E^*)$ such that for each
Borel subset $B$ of $X$
$$
\hat G_{y^*}(B)=y^*(\hat G(B))=y^*(G(B)).
$$
Hence since $\hat G=h\cdot |\hat G|_{\text {int }}$, it follows that for
each Borel subset $B$ of $X$
$$
\hat G_{y^*}(B)=\int_By^* h(x)\ d|\hat G|_{\text {int }}(x).
$$
This implies that for each $\varphi \in C(X,E)$
$$
\langle \varphi, \hat G_{y^*}\rangle=\int_Xy^*\langle h(x), \varphi
(x)\rangle d|\hat G|_{\text {int }}(x).\tag ****
$$
It follows that if $\sum\limits^n_{i=1}\varphi_i \otimes y^*_i$ is in
$C(X,E)\otimes_\epsilon F^*$ with $\parallel
\sum\limits^n_{i=1}\varphi_i \otimes y^*_i\parallel_\epsilon \leq 1$,
then for each $x$ in $X$
$$
|\langle h(x), \sum\limits^n_{i=1} \varphi_i(x)\otimes y^*_i\rangle
|\leq \parallel h(x)\parallel_{\text {int }}
$$
since
$$
\parallel \sum\limits^n_{i=1}\varphi_i(x)\otimes y_i^*\parallel_\epsilon
\leq \parallel \sum\limits^n_{i=1}\varphi_i\otimes
y^*_i\parallel_\epsilon
\leq
1
$$
It follows then from $(^{***})$ and $(^{****})$ that
$$
\parallel \sum\limits^n_{i=1}\langle T\varphi_i, y_i^*\rangle \parallel
\leq |\hat G|_{\text {int }}(X)<\infty
$$
Hence $J\circ T$ is an element of $I(C(X,E),F^{**})$ which implies that
$T\in I\left(C(X,E), F\right)$.  This completes the proof.
\bigskip
\noindent {\bf  Applications}
The study of integral operators on $C(X,E)$ spaces was motivated by some
problems that arose in [8] concerning nuclear operators on $C(X,E)$
spaces.  Recall that an operator $T$ between two Banach spaces $Y$ and
$Z$ is said to be a {\bf nuclear operator} if there exist sequences
$(y_n^*)$ and $(z_n)$ in $Y^*$ and $Z$ respectively, such that
$\sum\limits_n \parallel y_n^*\parallel \parallel z_n\parallel<\infty$
and for each $y\in Y$
$$
T_{\dsize y}=\sum\limits_n y_n^* (y) z_n.
$$
the nuclear norm is defined by
$$\parallel T\parallel_{nuc}=\text { inf }\left\{\sum\limits_n \parallel
y_n^* \parallel \parallel z_n\parallel \right\}
$$
where the infimum is taken over all sequences $(y_n^*)$ and $(z_n)$ such
that $T y=\sum\limits_n y_n^*(y)z_n$ for all $y\in Y$.  We shall denote by
$N(Y,Z)$  the space of nuclear operators from $Y$ to $Z$
under the nuclear
norm.  The study of nuclear operators on $C(X,E)$ spaces was initiated in
[1] where some of the known results in the scalar case were extended.
Nuclear operators on $C(X,E)$ spaces were also considered in [2] where it
was
shown that if $T:\ \ C(X,E)\longrightarrow F$ is a nuclear operator then
for each Borel subset $B$ of $X$, the operator $G(B):\ \ E\longrightarrow
F$ is nuclear.  In [8] it was shown that the representing measure $G$ of
a nuclear operator is in fact of bounded variation as a vector measure
taking its values in $N(E,F)$.  Hence all the above known results can
be summarized as follows:
\bigskip
\theorem { Proposition 4} { Let $X$ be a compact Hausdorff space and
let $E$ and $F$ be two Banach spaces.  If $T$ is a nuclear operator from
$C(X,E)$ into $F$ with representing measure $G$, then
\item {(i)}\ For any Borel subset $B$ of $X$ the operator
$G(B):E\rightarrow F$ is nuclear, and
\item {(ii)}\ $G$ is of bounded variation as a vector measure taking
its values in $N(E,F)$ under the nuclear norm.}
\medskip
Easy examples show that conditions (i) and (ii) above do not characterize
nuclear operators on $C(X,E)$ spaces.  As a matter of fact, it was shown
in [8] that counterexamples can be given as soon as the range space
$F$ fails to have the so-called Radon-Nikodym property (RNP) see [4].
The interesting question that arises is then:
\bigskip
\noindent {\bf Question 5}\ \ Assume that $F$ has the (RNP) and that $T:\
\ C(X,E)\longrightarrow F$ is an operator whose representing measure $G$
satisfies conditions (i) and (ii) of Proposition 4, is $T$ a nuclear
operator?
\bigskip
In what follows, we shall show how Theorem 3 can be used to give a
positive answer to Question 5 in case $F$ is assumed to be complemented
in its bidual.
\bigskip
\theorem { Theorem 6}{ Let $F$ be a Banach space with the
Radon-Nikodym property and assume that $F$ is complemented in its bidual.
A bounded linear operator $T:\ \ C(X,E)\longrightarrow F$ is nuclear
whenever its representing vector measure $G$ satisfies the following
conditions
\item {(i)} For each Borel subset $B$ of $X$, the operator $G(B):\ \
E\longrightarrow F$ is nuclear, and
\item {(ii)} $G$ is of bounded variation as a vector measure taking
its values in $N(E,F)$.}
\bigskip
\noindent {\bf Proof:}\ \ First note that for an arbitrary Banach space
$F$ any element $U$ in $N(E,F)$ is in $I(E,F)$ with $\parallel
U\parallel_{\text {int}}\leq \parallel U\parallel_{\text {nuc}}$.  Hence
by Theorem 3 an operator $T:\ \ C(X,E)\longrightarrow F$ that satisfies
conditions (i) and (ii) above is indeed integral.  Therefore when in
addition we assume that $F$ is complemented in its bidual $F^{**}$, the
operator $T$ becomes Pietsch integral (also called strictly integral) see
[4, p.235].  Under these conditions and if $F$ is assumed to have the
Radon-Nikodym property, it is well known [4, p.175] that $T$ becomes
a nuclear operator.  This completes the proof.
\bigskip
\noindent {\bf Remark}\ \ As we have just seen, Theorem 6 is a direct
application of the main result of this paper, it should be noted that
Theorem 6 gives a positive answer to Question 5 when $F$ is a dual space
and when $F$ is a Banach lattice.

\vfill\eject
\centerline {\bf Bibliography}
\bigskip
\item {[1]} G. Alexander, {\it Linear Operators on the Space of
Vector-Valued Continuous Functions}, Ph.D. thesis, New Mexico State
University, Las Cruces, New Mexico, 1972.

\item {[2]} R. Bilyeu and P. Lewis, {\it Some mapping properties of
representing measures}, Ann. Math. Pure Appl. CIX (1976), 273--287.

\item {[3]} G. Choquet, {\it Lecture on Analysis}, Vol. {\bf II},
Benjamin (1969).

\item {[4]} J. Diestel and J.J. Uhl, Jr., {\it Vector measures} Math.
Surveys {\bf 15} American Mathematical Society, (1977).

\item {[5]} N. Dinculeanu, {\it Vector measures}, Pergamon Press,
(1967).

\item {[6]} A. Grothendieck, {\it Produit Tensoriels topologiques et
Espaces Nucl\'eaires}.  Memoirs of the American Mathematical Society,
Vol. {\bf 16}, (1955).

\item {[7]} P. Saab, {\it Integral representation by boundary vector
measures}.  Canadian Math. Bull. {\bf 25 (2)} (1982), 164--168.

\item {[8]} B. Smith, {\it Some Bounded Linear Operators on the Spaces
$C(\Omega, E)$ and $A(K,E)$}, Ph.D. thesis, University of Missouri,
Columbia, Missouri, 1989.

\item {[9]} C. Swartz, {\it Absolutely summing and dominated operators
on spaces of vector valued continuous functions},Trans. Amer. Math. Soc.,
{\bf 179}, (1973), 123-132.

\vskip 1truein
\noindent {\bf University of Missouri}

\noindent {\bf Department of Mathematics}

\noindent {\bf Columbia, MO  65211}

\vfill\eject\end